\documentclass[11pt]{article}
\pdfoutput=1
\usepackage{amsmath}
\usepackage{amssymb}
\usepackage{amsthm}
\usepackage{listings}
\usepackage{fancyhdr}
\usepackage[backend=bibtex, style=numeric, sorting=none]{biblatex}
\usepackage{multicol}
\newtheorem{theorem}{Theorem}[section]
\newtheorem{definition}[theorem]{Definition}
\newtheorem{example}[theorem]{Example}

\begin{document}
	\title{An iterative algorithm for approximating roots of integers}
	\author{William Gerst}
	\date{January 8, 2021}
	\maketitle

\begin{abstract}
	We explore an algorithm for approximating roots of integers, discuss its motivation and derivation, and analyze its convergence rates with varying parameters and inputs. We also perform comparisons with established methods for approximating square roots and other rational powers.
\end{abstract}

\section{Introduction}
	Fixed-point iteration is a technique commonly used for approximating irrational constants, such as calculating the golden ratio $\phi\approx1.618$ or finding square roots via the Babylonian method. More generally, the Newton-Raphson method is an iterative process used for estimating the roots of functions, which can be applied to find the roots of certain integers. For example, searching for the real root of $f(x)=x^3-5$ will yield $x=\sqrt[3]5$.
	
	In this paper, we will introduce an algorithm for estimating rational powers of integers that employs fixed-point iteration, provide examples of its use, and compare it to the existing methods mentioned above.
	
	It is important to note the definition of a ``fixed point" in the context of functions and iterative processes, as well as what it means for a fixed point to be ``stable."
	
\begin{definition}
	A \textbf{fixed point} $x_0$ of a function $f$ is a value in the domain of $f$ such that $f(x_0)=x_0$.
\end{definition}

	In other words, a fixed point of a function occurs wherever the input and output values are equal.

\begin{example}
	The values $\displaystyle\phi=\frac{1+\sqrt5}{2}=1.618\ldots$ and $\displaystyle\Phi=-\frac1\phi\approx-0.618$ are the two fixed points of the function $\displaystyle f(x)=1+\frac1x$ because $f(\phi)=\phi$ and $f(\Phi)=\Phi$.
\end{example}

	Not all fixed points are of equal importance, however. One important distinction to be made for a given fixed point is whether it can be classified as ``stable" or ``unstable," which describes how the function acts in the neighborhood of this point.

\begin{definition}
	A fixed point $c$ of a function $f$ is \textbf{stable} if and only if the sequence $x_{n+1}=f(x_n)$ converges to $c$ with $x_1=c+\varepsilon$ for arbitrarily small $\varepsilon$. Conversely, $c$ is \textbf{unstable} if and only if the sequence $x_{n+1}=f(x_n)$ diverges from $c$ with small initial offset $\varepsilon$.
\end{definition}

\renewcommand{\arraystretch}{1.2}
	The stability of fixed points of a function $f(x)$ can be identified by the magnitude of the derivative at those points \cite{stability}.
\begin{center}
\begin{tabular}{l|l}
	unstable & $\left\lvert f'(x)\right\rvert>1$ \\
	neutral & $\left\lvert f'(x)\right\rvert=1$ \\
	stable & $\left\lvert f'(x)\right\rvert<1$ \\
	superstable & $\left\lvert f'(x)\right\rvert=0$
\end{tabular}
\end{center}

	When a fixed point is ``superstable," it is said that the algorithm ``converges (at least) quadratically" because the error in the approximation $x_{n+1}$ must be some quadratic function of the error in $x_n$.
	
\section{Motivation}

	In the Newton-Raphson method, we see an example of how fixed-point iteration can be used to narrow in on an estimated value for the root of a given integer. Suppose that we begin the process in a slightly different way.
	
	For clarity of this example, suppose we are trying to find an estimated value of $\sqrt 7$. Rather than defining the function $f(x)=x^2-7$ and trying to find the positive root, we simply consider the equation $r^2=7$ for some $r$. Ideally, we could divide each side of the equation by $r$ to get $r=7/r$ and iterate as $r_{n+1}=7/r_n$ with arbitrary $r_1$ to get closer and closer to $r_n=\sqrt 7$ as $n\to\infty$. There is, however, a fatal flaw with this approach: the fixed point of $r=\sqrt7$ has only neutral stability, so the value does not converge as we might hope. Instead, it can be easily seen that any nonzero starting value of $r_1$ will simply yield values alternating between $r_1$ and $7/r_1$.
	
	This flaw can be remedied if we offset $r$ by some constant value. Suppose that we rewrite $r=2-c$. (As for where the 2 and the negative sign come from, that will be explained shortly.) Then, we have the initial equation $(2-c)^2=4-4c+c^2=4+c(c-4)=7\iff c=\frac{7-4}{c-4}$ and can turn this into the fixed-point iteration with
	$$c_{n+1}=\frac{3}{c_n-4}.$$
	
	Taking the derivative of the right-hand side with respect to $c$ at the fixed point of $c=2-\sqrt7$ gives
	$$\frac{\mathrm d}{\mathrm dc}\left(\frac{3}{c-4}\right)=\frac{-3}{(c-4)^2}=\frac{-3}{(-2-\sqrt7)^2}=\frac{-3}{4+4\sqrt7+7}\approx-0.139.$$
	
	The absolute value of the derivative at the fixed point is less than 1, so it is stable! This means that repeatedly iterating on the value of $c$ will cause the value to tend towards $2-\sqrt7$, at which point we may subtract $2-c$ to recover our approximation for $\sqrt7$ as desired.
	
	This is where the earlier rewriting of $r=2-c$ is important. Had we chosen a different constant or sign, it is entirely possible that the fixed point would not have been stable, ruining the results. In the following section, we generalize our results and make clear exactly how to choose this constant to make the algorithm work (and optimize its rate of convergence).
	
\section{Derivation of general algorithm}

	First, we want to generalize to arbitrary input $x\in\mathbb R^+$, not just the number 7 used in the previous example. (Note that each of the examples given in this paper use integer $x$, but there is no limitation that prevents calculating $\sqrt\pi$, for example. The main drawback is simply that it is more difficult to calculate by hand and of less interest to the general reader.) Start with the quadratic equation
	\begin{equation}
		\left(\frac b2-c\right)^2=x
	\end{equation}
	where $b/2$ is our ``offset" value and we can solve to get $\frac b2-c=\pm\sqrt x$. Next, rewrite equation (1) in expanded form as
	\begin{equation*}
		c^2-bc+\frac{b^2}4=x
	\end{equation*}
	and perform algebraic manipulation to obtain
	\begin{equation*}
		c=\frac{x-(b^2/4)}{c-b}
	\end{equation*}
	which can be easily transformed into the fixed-point iteration of
	\begin{equation}
		c_{n+1}=\frac{x-(b^2/4)}{c_n-b}.
	\end{equation}
	
	From here, it is easily verified that the two fixed points for $c$ are given by $c_-=\frac b2-\sqrt x$ and $c_+=\frac b2+\sqrt x$. Determining which of these is stable is fundamental in knowing how to proceed.
	
	We define $f(c)=\frac{x-(b^2/4)}{c-b}$ and take the absolute value of the derivative to get
	$$|f'(c)|=\left|\frac{x-(b^2/4)}{(c-b)^2}\right|.$$
	
	We first check the fixed point of $c_+=\frac b2+\sqrt x$. Substituting yields
	$$|f'(c_+)|=\left|\frac{x-(b^2/4)}{(\sqrt x-b/2)^2}\right|=\left|\frac{x-(b^2/4)}{x-b\sqrt x+(b^2/4)}\right|.$$
	
	Consider the case in which we have picked $b$ such that $0<b<2\sqrt x$. We can rewrite $\sqrt x=b/2+\varepsilon$ for some $\varepsilon>0$. Then, we can simplify to
	$$|f'(c_+)|=\left|\frac{b\varepsilon+\varepsilon^2}{\varepsilon^2}\right|>1.$$
	
	In this case, then, the fixed point $c_+$ is unstable. Consider instead the case of $0<2\sqrt x<b$. We rewrite $\sqrt x=b/2-\varepsilon$ for some $0<\varepsilon<b/2$. Then, we can simplify to
	$$|f'(c_+)|=\left|\frac{-b\varepsilon+\varepsilon^2}{\varepsilon^2}\right|=\left|1-\frac{b}{\varepsilon}\right|>1,$$
	where the final inequality holds true because of the bounds on $\varepsilon$. This is unstable once more.
	
	Next, let us consider the case of $-2\sqrt x<b<0$. Rewriting $\sqrt x=-b/2+\varepsilon$ for some value $\varepsilon>0$ gives
	$$|f'(c_+)|=\left|\frac{-b\varepsilon+\varepsilon^2}{(-b+\varepsilon)^2}\right|=\left|\frac{\varepsilon}{\varepsilon-b}\right|<1.$$
	Unlike the first two cases, this shows that the fixed point $c_+$ is stable for such a selection of $b$.
	
	Fourth, we consider the case of $b<-2\sqrt x<0$. Rewriting $\sqrt x=-b/2-\varepsilon$ for some value $0<\varepsilon<-b/2$ gives
	$$|f'(c_+)|=\left|\frac{b\varepsilon+\varepsilon^2}{(-b-\varepsilon)^2}\right|=\left|\frac{\varepsilon}{\varepsilon+b}\right|<1.$$
	Once again, the fixed point is stable.
	
	The last inequality here may not be immediately clear, but it can be observed that the value of the denominator lies within the open interval of $b$ (which is greater in magnitude than any possible value of $\varepsilon$) to $b/2$, at which point the numerator and denominator are equal in magnitude.
	
	Now, we address the few remaining special cases. When we have the exact equality $b=2\sqrt x$, we see $|f'(c_+)|\to\infty$, making it unstable. When $b=0$, we have that $|f'(c_+)|=1$, indicating a neutral stability (the same situation described as a ``flaw" in the beginning of section 2). Finally, when we have $b=-2\sqrt x$, we see $|f'(c_+)|=0$, making the point superstable.
	
	Nearly identical calculations reveal that whenever the point $c_+$ is stable, the point $c_-$ is unstable, and vice versa. A summary table of all the stability results above for the different possible values of $b$ is provided below:
	
	\begin{center}
		\begin{tabular}{c|c|c} 
			& $c_-=b/2-\sqrt x$ & $c_+=b/2+\sqrt x$ \\
			\hline
			$b<-2\sqrt x$ & unstable & stable \\
			$b=-2\sqrt x$ & unstable & superstable \\ 
			$-2\sqrt x<b<0$ & unstable & stable \\
			$b=0$ & neutral & neutral \\ 
			$0<b<2\sqrt x$ & stable & unstable \\
			$b=2\sqrt x$ & superstable & unstable \\ 
			$2\sqrt x<b$ & stable & unstable
		\end{tabular}
	\end{center}

	Essentially, the key takeaway is that there is a kind of symmetry about $b=0$ with respect to the stability of these fixed points. To make things simpler, we can simply choose to consider strictly positive $b$ and deal only with $c_-$, or we may choose to consider only negative $b$ and the corresponding $c_+$; either is equally valid. In this analysis, we opt with the former.
	
	\subsection{How to put the square root algorithm into practice}
	
	A quick glance at the summary table provides instructions on how to select the constant $b$: choose a positive value as close to $2\sqrt x$ as possible to maximize the convergence rate. (Estimating $\sqrt x$ is what we seek to do with this algorithm, so coming up with a good initial guess is not necessarily an easy task. It can be noted, however, that for an input integer $x$ with $n$ digits, the square root of $x$ will have approximately $n/2$ digits, which can be used to rapidly generate a ballpark figure.)
	
	Once we have chosen an appropriate value of $b$, we want to choose some initial value $c_1$ for our iteration. Fortunately, the convergence does not depend on this value, except for that the denominator of (2) must never be zero. To avoid this, set $c_1$ to some number not equal to $b$, and if by chance a division by zero does occur, simply restart the process with a new value of $c_1$. (Ideally, $c_1\approx b/2-\sqrt x$ will provide the best accuracy in the fewest iterations. We can use the same trick as before to estimate $\sqrt x$ here if desired, or simply note that this estimate gives $b/2\approx\sqrt x$ and set $c_1=0$.)
	
	After values of $b$ and $c_1$ are both chosen, iterate with equation (2) until achieving the desired level of accuracy (in terms of unchanging decimal places, etc.).
	
	To finish the process, calculate $b/2-c_n$ to get an estimate for $\sqrt x$.
	
\subsection{Example of square-root application}
	
	Suppose we wish to calculate an approximation of $\sqrt5$. Since $4\approx5$, it stands to reason that $\sqrt4=2\approx\sqrt5$, so we will choose the value $b=4$ as our guess. Since we have a fairly good estimate for $b$, it would make sense to choose $c_1=0$. For the purpose of demonstrating the convergence speed even in non-optimal cases, however, we will choose $c_1=10$ here instead. Now, we use the recursive formula
	\[c_{n+1}=\frac{x-4}{c_n-4}\]
	to calculate $c_n$ up to our desired level of precision. Afterwards, we will evaluate $2-c_n$ to obtain an approximation of $\sqrt 5$. The results of these calculations are provided below:
	
	\begin{center}	
		$\begin{array}{c|c}
			n & c_n \\
			\hline
			1 & 10 \\
			2 & 0.16666666\dots \\
			3 & -0.26068956\dots \\
			4 & -0.23469387\dots \\
			5 & -0.23614457\dots \\
			6 & -0.23606370\dots \\
			7 & -0.23606821\dots \\
			8 & -0.23606796\dots \\
			9 & -0.23606797\dots \\
			10 & -0.23606797\dots \\
		\end{array}$
	\end{center}
	
	Finally, we conclude by computing $2-c_{10}\approx2.23606797745$, to be compared against the desired value of $\sqrt5\approx2.23606797750$.
	
	Interestingly, this specific process produces the continued fraction
	\[-[0;\bar4]=2-\sqrt5=\cfrac1{-4+\cfrac1{-4+\cfrac1{-4+\cfrac1\cdots}}}\]
	that clearly demonstrates why $2-c\approx2-(2-\sqrt5)=\sqrt5$.
	
	\subsection{Attempting to expand to $n^{th}$ roots}
	
	It is clear how to extend this reasoning to more than simply square roots. To estimate a power $\sqrt[d]x$, we start with the analogous form of equation (1):
	\begin{equation}
		\left(\frac b2-c\right)^d=x\iff\frac b2-c=\sqrt[d]x.
	\end{equation}

	To make analysis easier, we can swap the sign of $c$ to give
	$$\left(\frac b2+c\right)^d=x\iff\frac b2+c=\sqrt[d]x.$$
	
	Binomial expansion of the left equality yields
	$$x=\left(\frac b2\right)^d+{d\choose 1}\left(\frac b2\right)^{d-1}c+{d\choose 2}\left(\frac b2\right)^{d-2}c^2+\cdots+{d\choose d-1}\left(\frac b2\right)c^{d-1}+c^d.$$
	
	All terms, except for the first, in the binomial expansion have at least one factor $c$. Then, using the same algebraic manipulations as before, it is trivial to find that
	$$x-(b/2)^d=c\cdot\left(c^{d-1}+\cdots\right)=c\cdot\left(\frac{(b/2+c)^d-(b/2)^d}c\right).$$
	
	Cross-multiplying, it follows that the generalized form of (2) is given by
	\begin{equation}
		c_{n+1}=\frac{c_n\cdot\left(x-\left(\frac b2\right)^d\right)}{\left(\frac b2+c_n\right)^d-\left(\frac b2\right)^d}.
	\end{equation}
	where the fixed point for $c$ satisfies
	$$\sqrt[d]x\approx\frac b2+c.$$
	
	The equations (1) and (2) can be readily recovered by setting $d=2$ and replacing $c$ with $-c$ in formulas (3) and (4) and simplifying.
	
	Setting $d=3$, on the other hand, results in the equation
	$$c_{n+1}=\frac{x-(b^3/8)}{c_n^2+3bc_n/2+3b^2/4}.$$
	Note that we can choose to avoid doing this simplification and just use equation (4) directly, which is especially helpful for large $d$.
	
	However, trying this equation with test values may reveal a disappointing truth: the equation is not as versatile as the square root special case. As an exercise, attempt to calculate $\sqrt[3]{1250}$ with the values $(b,c_1)=(7,1)$ or $(b,c_1)=(7,10)$. They don't converge to the correct value! On the other hand, testing with $(b,c_1)=(70,1)$ or $(b,c_1)=(70,10)$ does produce the correct value of about $10.772$. To understand why this occurs, we must return to the idea of fixed point stability.
	
	In the quadratic special case, we found the fixed points $c$ of the algorithm as dictated by the $b$ and $x$ inputs. These were of the form $c_-=b/2-\sqrt x$ and $c_+=b/2+\sqrt x$. Now, in the higher-degree case, we have performed a sign swap, so these fixed points are of the form $c=\sqrt[d]x-b/2$ (with a $\pm$ on the radical term in the case of even $d$). Then, to analyze the fixed point stability, we differentiate the right-hand side of (4) with respect to $c$ at the fixed point(s). Finding where the magnitude of the derivative crosses one, we will see for which values of $b$ the fixed point(s) switch from being stable to unstable. In the example with $\sqrt[3]{1250}$, this yields values of $b\approx7.88578$ and $b\approx-29.43013$, where the algorithm fails to converge to the correct value within the interval but succeeds outside the interval. (Solving the resulting quadratic equation in $b$ analytically gives $b=\sqrt[3]x\cdot(\pm\sqrt3-1)$, which supports the aforementioned empirical values.)
	
	The cases of even higher degree work even less frequently. In computer-based testing, it seems that the algorithm works but only for specific values of $b$ inside tight intervals. This is a stark contrast with the quadratic case, which required only that $b$ have a particular sign.
	
	In the seventh-degree case, for example, the least positive value of $b$ for which calculating $\sqrt[7]{1250}$ seems to work is approximately $4.1755$, while the lowest possible positive value of $b$ for estimating $\sqrt[7]{12500}$ is approximately $5.8019$. Guessing that the seventh powers of these bounds are related to their corresponding values of $x$ (analogous to the cubic example), we see that the ratio $b_{lower}^7/x\approx17.704$ in each case. Extrapolating to $x=125000$ predicts $b_{lower}\approx8.0617$, which can be verified via computer. Evidently, this trend is real, although the exact mechanism by which the number 17.704 can be generated remains elusive.
	
	Unfortunately, there does not seem to be a clear path to a general analytical solution to this problem of determining the interval(s) of validity. It should be noted that whereas there was no upper bound on the valid range of $b$ values for the cubic case, the same does not seem to hold true for higher-degree examples. Fortunately, computer testing suggests that even at high degrees, there are at least some values of $b$ for which this process works, so a solution could recover some of the usefulness of this technique.

\section{Comparison to existing methods}
	There are two connections we will draw. Namely, we will discuss connections to continued fractions and the Newton-Raphson method mentioned before.
	
	We recall that the special case $d=2$ in equation (4) results in (2):
	\[c_{n+1}=\frac{x-(b^2/4)}{c_n-b}\]
	In this case, we see that the recurrence relation allows for an infinite generalized continued fraction to be constructed \cite{contfrac}.
	\[c=\frac b2-\sqrt x=\frac{x-(b^2/4)}{-b+\cfrac{x-(b^2/4)}{-b+\cfrac{x-(b^2/4)}{\ddots}}}\]
	In 3.2, the numerators were each 1, making a ``simple" continued fraction.

	Next, the Newton-Raphson method is another technique used for approximating roots of positive real numbers, and it also relies on fixed-point iteration in its calculations.
	
	Suppose we want to find an approximation for $\sqrt7$. We know that we can achieve this by searching for the positive root of $f(x)=x^2-7$. By setting $x_1=1$ and performing successive linear approximations with
	\begin{equation}
		x_{n+1}=x_n-\frac{f(x_n)}{f'(x_n)}=x_n-\frac{x_n^2-7}{2x_n},
	\end{equation}
	we will have $x_n$ converge to the desired value. (This is equivalent to the Babylonian method, after algebraic simplification.)
	
	In fact, we see that at the fixed point $x_n=\sqrt7$, the derivative of the right-hand side with respect to $x_n$ is zero. This indicates that the Newton-Raphson method converges at least quadratically (the number of accurate decimal places roughly doubles in each iteration), which is better than the method explored in this paper, which only converges linearly (since the derivative is generally nonzero, as discussed in the section on stability). This can be explained by the theory of Taylor expansion \cite{quadconv}.
	
	Thus, the Newton-Raphson method and higher-order analogues like Halley's method experience a much faster rate of convergence, with the trade-off being that elementary differentiation of the polynomial $f(x)$ is necessary for their computation. Since this polynomial is always of the form $f(x_n)=x_n^d-k$ with degree $d$ and constant $k$, we know that $f'(x_n)$ in right-hand side of (5) will always be equal to $dx_n^{d-1}$ by the power rule. For Halley's method and higher-order root-finding algorithms, the higher-order derivatives can be similarly computed without much effort.
	
\printbibliography

\end{document}